\documentclass[11pt]{article}
\usepackage[centertags]{amsmath}
\usepackage{amsfonts}
\usepackage{amssymb}
\usepackage{amsthm}
\usepackage{amsmath}
\usepackage{amsxtra}
\newcommand{\BN}{\Bbb{N}}
\newcommand{\BZ}{\Bbb{Z}}

\newcommand{\D}{\displaystyle}

\newcommand{\fa}{\frak{a}}

\begin{document}
\pagestyle{plain}
\thispagestyle{empty}
 \title
{Finiteness of Composed Local Cohomology modules }
\author{Fatemeh Dehghani-Zadeh}
\date{}

\maketitle \footnote{{\em keywords:}
Generalized local cohomology; Finiteness; Serre subcategories. }

\noindent {\bf\large\hspace{0 cm} \noindent Abstract. ~}
 Cofiniteness of the generalized local cohomology modules $H^{i}_{\mathfrak{a}}(M,N)$ of two $R$-modules $M$ and $N$ with respect to an ideal $\mathfrak{a}$ is studied for some $i^{,}s$ with a specified property.
Furthermore, Artinianness of
$H^{j}_{\mathfrak{b}_{0}}(H_{\mathfrak{a}}^{i}(M,N))$ is investigated by using
the above result, in certain graded situations, where $\mathfrak{b}_{0}$ is an ideal of $R_{0}$ such that $\mathfrak{b}_{0}+\mathfrak{a}_{0}$ is an $\mathfrak{m}_{0}$-primary ideal.
 \section*{1. Introduction}
 \indent Generalized local cohomology was given in the local case by
 J. Herzog [7] and in the more general case by Bijan-Zadeh [2].
 Let $R$ be a commutative Noetherian ring ( not necessarily local) with identity, $\mathfrak{a}$ an ideal of $R$ and let $M, N$ be two $R$-modules. For an integer $i\geq0$, the $i$th generalized local cohomology module $H^{i}_{\mathfrak{a}}(M,N)$
 is defined by $\D{H_{\mathfrak{a}}^{i}(M,N)
 =\varinjlim_{n\in{\BN}}\hbox{Ext}_{R}^{i}({M}{/{\mathfrak{a}}^nM},N)}$. With $M=R$,
 we obtain the ordinary local cohomology module
 $H^{i}_{\mathfrak{a}}(N)$ of $N$ with respect to $\mathfrak{a}$ which was
 introduced by Grothendieck.\\
Recall that a class $S$ of $R$-modules is a Serre subcategory of
the category of $R$-modules, when it is closed under taking
submodules, quotients and extensions. One can see that the subcategories of minimax and $\frak{a}$-cofinite modules, weakly laskerian $R$-modules, $R$-modules with finite support are examples of Serre classes.\\This paper is divided into four sections. In the second
section of the paper, we study some results on Serre classes by using spectral sequences. In the section 3, we investigate the Cofiniteness and minimaxness property of generalized local cohomology modules.\\ Bahmanpour and Naghipour in [1, Theorem 2.6] showed that if $t$ is a non-negative such that $\hbox{Supp}H^{i}_{\frak{a}}(N)\leq1$ for $i<t$ then $H^{i}_{\fa}(N)$ is $\frak{a}$-cofinite for all $i<t$ and cuong, Goto and Hoang proved that if  $\hbox{Supp}H^{i}_{\frak{a}}(M, N)\leq 1$ for all $i<t$ then $H^{i}_{\fa}(M, N)$ is $\frak{a}$-cofinite for all $i< t$ [5, Theorem 1.2].
Our main aim in this section is to establish the following theorem:\\
Theorem 1.1. Let $t$ be a non-negative integer such that $\hbox{Supp}H^{i}_{\frak{a}}(M, N)$ is finite set for all $i<t$. Then $H_{\mathfrak{a}}^{i}(M, N)$ is $\frak{a}$-cofinite and minimax $R$-module for all $i<t$.\\
Throughout section 4, $R=\bigoplus_{n\geq0}R_{n}$ is a graded
commutative Noetherian ring, where the base ring $R_{0}$ is a
commutative Noetherian local ring with maximal ideal
$\mathfrak{m}_{0}$. Moreover, we use $\mathfrak{a}_{0}$ to denote
a proper ideal of $R_{0}$ and we set $R_{+}=\bigoplus_{n>0}R_{n}$
the irrelevant ideal of $R$,
$\mathfrak{a}=\mathfrak{a}_{0}+R_{+}$, and
$\mathfrak{m}=\mathfrak{m}_{0}+R_{+}$ and $\frak{b}_{0}$ is an ideal of $R_{0}$ such that $ \frak{a}_{0}+\frak{b}_{0}$ is $\frak{m}_{0}$-primary ideal. Also, we use
$M=\bigoplus_{n\in {\BZ}}M_{n}$ and $N=\bigoplus_{n\in
{\BZ}}N_{n}$ to denote non-zero, finitely generated graded
$R$-modules. It is well known that, the $i$th generalized local
cohomology module $H_{\mathfrak{a}}^{i}(M,N)$ inherits natural
grading for each $i\in \BN_{0}$ (where $\BN_{0}$ denotes the set
of all non-negative integers). In section 4, using the results of section 3, we study the Artinianness and cofiniteness  of $R$-modules $H^{j}_{\frak{b}_{0}}(H^{i}_{\frak{a}}(M, N))$. More precisely we shall show that:\\
Theorem 1.2. Let $T$ be an $\frak{a}$-torsion graded minimax $R$-module and $ \frak{a}_{0}+\frak{b}_{0}$ is $\frak{m}_{0}$-primary ideal. Then $H^{i}_{\frak{b}_{0}}(T)$ and $Tor_{i}^{R_{0}}(\frac{R_{0}}{\frak{b}_{0}}, T)$ are Artinian for all $i\geq 0$.\\
Theorem 1.3. Let $T$ be an $\frak{a}$-torsion and $\frak{a}$-cofinite $R$-module and $ \frak{a}_{0}+\frak{b}_{0}$ is $\frak{m}_{0}$-primary ideal. Then $H^{i}_{\frak{b}_{0}}(T)$ is Artinian and $\frak{a}$-cofinite and for all $i\geq 0$.\\
Theorem 1.4. Let $t$ be a non-negative integer such that $\hbox{Supp}H^{i}_{\frak{a}}(M, N)$ is finite set for all $i<t$ and $ \frak{a}_{0}+\frak{b}_{0}$ is $\frak{m}_{0}$-primary ideal. Then $H^{j}_{\frak{b}_{0}}(H^{i}_{\frak{a}}(M, N))$  is Artinian and $\frak{a}$-cofinite and for all $j\geq 0$ and $i<t$.\\
Theorem 1.5. Let $\hbox{gdepth}(M/\frak{a}M , N)= \hbox{cd}(M, N)$ and $ \frak{a}_{0}+\frak{b}_{0}$ is $\frak{m}_{0}$-primary ideal. Then $H^{j}_{\frak{b}_{0}}(H^{i}_{\frak{a}}(M, N))$  is Artinian and cofinite and for all $j\geq 0$ and $i\geq 0$.\\
Throughout this paper, $R$ will always be a commutative Noetherian ring. For a representable $R$-module $T$ we denote by $\hbox{Att}_{R}(T)$ the set of attached prime ideals of $T$. We shall use $\hbox{Max}(R)$ to denote the set of all maximal ideals of $R$. \\For any unexplained notation and terminology we refer the reader to [3], [4] and [15].\\
 \section*{2. Serre classes and Spectral sequences}
 This section is started with the following theorem.\\
{\bf Theorem 2.1}. (cf [15, Theorem 10.47]). Let $ A\overset{G}{\longrightarrow} B\overset{F}\longrightarrow C$ be covariant additive functors, where $A, B$, and $C$ are abelian categories with enough injectives. Assume that $F$ is left exact and that $GE$ is right $F$-acyclic for every injective object $E$ in $A$. Then, for every object $T$ in $A$, there is a third quadrant spectral sequence with $$E_{2}^{p,q}=(R^{p}F)(R^{q}G)T\underset{p}\Longrightarrow R^{n}(FG)T=E^{n}.$$
{\bf Theorem 2.2}. Let the situation be as in Theorem 2.1. If $(R^{p}F)(R^{q}G)T=0$ for all $p>r$ or $q>s$ or $p<c$ or $q<d$, then the following statements are proved:
\begin{itemize} \item[(i)] $(R^{r}F)(R^{s}G)T\cong R^{r+s}(FG)T,$
\item[(ii)] $(R^{c}F)(R^{d}G)T\cong R^{c+d}(FG)T.$
\end{itemize}
{\bf Proof}. (i) It follows from Theorem 2.1 the Grothendieck$^{,}$s spectral sequence $$E_{2}^{p,q}=(R^{p}F)(R^{q}G)T\underset{p}\Longrightarrow R^{p+q}(FG)T=E^{n}.$$ We consider the exact sequence $E_{k}^{r-k,s+k-1}\longrightarrow E_{k}^{r,s}\longrightarrow E_{k}^{r+k,s+1-k}.$
Since $E_{k}^{r-k,s+k-1}=E_{k}^{r+k,s-k+1}=0$ for all $k\geq 2$, we get $E_{2}^{r,s}=E_{3}^{r,s}=\cdots=E_{\infty}^{r,s}.$\\
On the other hand, there is a filtration $\varphi$ of $E^{r+s}$ with $$0=\varphi^{r+s+1}E^{r+s}\subseteq \varphi^{r+s}E^{r+s}\subseteq \cdots
\subseteq\varphi^{1}E^{r+s}\subseteq \varphi^{0}E^{r+s}=E^{r+s}$$ such that $E_{\infty}^{i,r+s-i}=\frac{\varphi^{i}E^{r+s}}{\varphi^{i+1}E^{r+s}}$ for all $0\leq i \leq r+s$. As $E_{2}^{i,r+s-i}=0$ for all $i\neq r$, we have $\varphi^{r+1}E^{r+s}=\varphi^{r+2}E^{r+s}=\cdots =\varphi^{r+s+1}E^{r+s}=0$ and $\varphi^{r}E^{r+s}=\varphi^{r-1}E^{r+s}=\cdots =\varphi^{0}E^{r+s}=E^{r+s}$. It follows $E_{\infty}^{r,s}=\frac{\varphi^{r}E^{r+s}}{\varphi^{r+1}E^{r+s}}\cong \varphi^{r}E^{r+s}\cong E^{r+s}.$ So $(R^{r}F)(R^{s}G)T\cong R^{r+s}(FG)T.$\\
(ii) By an argument similar to the $(i)$ we get $E^{c,d}_{2}=E^{c,d}=\cdots=E^{c,d}_{\infty}.$ We have a filtration $\varphi$ of $E^{c+d}=R^{c+d}(FG)T$ with $0=\varphi^{d+c+1}E^{c+d}\subseteq \varphi^{d+c}E^{d+c}\subseteq \cdots \subseteq \varphi^{1}E^{c+d}\subseteq \varphi^{0}E^{c+d}=E^{c+d}$ such that $\varphi^{c+1}E^{d+c}=\varphi^{c+2}E^{d+c}=\cdots =\varphi^{c+d+1}E^{c+d}=0$ and $\varphi^{c}E^{c+d}=\varphi^{c-1}E^{c+d}=\cdots =\varphi^{0}E^{c+d}=E^{c+d}.$ It follows $E_{\infty}^{c,d}=\frac{\varphi^{c}E^{c+d}}{\varphi^{c+1}E^{c+d}}=\varphi^{c}E^{c+d}=E^{c+d}$. This proves the claim.

{\bf Theorem 2.3}. Let the situation be as in Theorem 2.1. If $S$ is Serre class and $(R^{p}F)(R^{q}G)T$ is in $S$ for all $q<r$, then $R^{n}(FG)T$ is in $S$ for $n<r$.\\
{\bf Proof}. By Theorem 2.1, there is a Grothendieck$^{,}$s spectral sequence $$E_{2}^{p,q}=(R^{p}F)(R^{q}G)T\underset{p}\Longrightarrow R^{p+q}(FG)T=E^{n}.$$ For all $i\geq 2$, $t<r$ and $p\geq 0$, we consider the exact sequence $$ 0\longrightarrow \hbox{ker} d_{i}^{p,t}\longrightarrow
E_{i}^{p,t}\longrightarrow E_{i}^{p+i,t-i+1}.  \;\;\;    (1)$$ Since
$E_{i}^{p,t}=\frac{\hbox{ker}d_{i-1}^{p,t}}{\hbox{im}
d_{i-1}^{p-i+1,t+i-2}}$ and $E_{i}^{p,j}=0$ for all $j<0$, we use $(1)$ to obtain $\hbox{ker}d_{t+2}^{i,t-i}\cong
E^{i,t-i}_{t+2}\cong\ldots\cong E_{\infty}^{i,t-i}$ for all $0\leq
i\leq t$ . In addition, there exists a finite filtration
$$0=\varphi^{t+1}E^{t}\subseteq \varphi^{t}E^{t}\subseteq\ldots
\subseteq\varphi^{1}E^{t}\subseteq
\varphi^{0}E^{t}=E^{t}$$such that
$E_{\infty}^{i,t-i}=\frac{\varphi^{i}E^{t}}{\varphi^{i+1}E^{t}}$
for all $0\leq i\leq t$.\\ Now, the exact sequence
$$0\longrightarrow \varphi^{i+1}E^{t}\longrightarrow
\varphi^{i}E^{t}\longrightarrow E_{\infty}^{i,t-i}\longrightarrow
0       \;\; \;\; ( 0\leq i\leq t)$$  in conjunction with
$E^{i,t-i}_{\infty}\cong \hbox{ker}d_{t+2}^{i,t-i}\subseteq
\hbox{ker}d_{2}^{i,t-i} \subseteq E_{2}^{i,t-i}$ yields
$E^{i}$ is in $S$ for all $0\leq i<r$.\\
 {\bf Theorem 2.4}. Let the situation be as in Theorem 2.1. If $S$ is Serre class and $(R^{p}F)(R^{q}G)T$ is in $S$ for all $q<r$ and $R^{i}(FG)T$ is in $S$ for all $i\geq 0$, then $(R^{p}F)(R^{r}G)T$  is in $S$ for all $p=0, 1$.\\
{\bf Proof}. By Theorem 2.1, there is a Grothendieck$^{,}$s spectral sequence $$E_{2}^{p,q}=(R^{p}F)(R^{q}G)T\underset{p}\Longrightarrow R^{p+q}(FG)T=E^{n}.$$ Also, there is a bounded filtration
$0=\varphi^{t+1}E^{t}\subseteq\varphi^{t}E^{t}\subseteq \cdots \subseteq \varphi^{1}E^{t}\subseteq\varphi^{0}E^{t}=E^{t}$ such that $E^{i,t-i}_{\infty}\cong\frac{\varphi^{i}E^{t}}{\varphi^{i+1}E^{t}}$ for all $0\leq i \leq t$ and hence $E^{p,q}_{\infty}$ is in $S$ for all $p, q$. Note that $E_{\infty}^{p,q}=E_{r}^{p,q}$ for large $r$ and each $p$ and $q$. It follows that there is an integer $\ell\geq 2$ such that $E_{r}^{p,q}$ is in $S$ for all $r\geq \ell$. We now argue by descending induction on $\ell$. Now, assume that $2<\ell <r$ and that the claim holds for $\ell$. Since $E^{p,q}_{r}$ is in a subquotient of $E^{p,q}_{2}$ for all $p, q \in\BN_{0}$, the hypotheses give $E^{p+r,t-r+1}_{r}$ is in $S$ for all $r\geq 2$. In addition, $E^{p,t}_{\ell}=\frac{\hbox{ker}d^{p,t}_{\ell-1}}{\hbox{im}d^{p-\ell+1}_{\ell-1}}$ and $\hbox{im}d^{p-\ell+1}_{\ell-1}=0$ for $p=0, 1$, it follows that $\hbox{ker}d^{p,t}_{\ell-1}$ is in $S$ for all $\ell>2$ and $p=0, 1$. Let $r\geq 2$ and $p=0, 1$, we consider the sequence $$0\longrightarrow \hbox{ker}d^{p,t}_{r} \longrightarrow E^{p,t}_{r} \longrightarrow E^{p+r,t-r+1}_{r}.$$ Since both $\hbox{ker}d^{p,t}_{\ell-1}$ and $E^{p+r,t-r+1}_{\ell-1}$ are in $S$, it follows that $E^{p, t}_{\ell-1}$ is in $S$ for $p=0, 1$. This completes the inductive step.
\section*{3. Minimax and Cofinite modules}
We keep the notation and hypotheses given in the introduction and continue with the following definition\\
{\bf Definition 3.1}. \begin{itemize} \item[(i)] An $R$-module $T$ is said to be $\frak{a}$ -cofinite if $\hbox{Supp}T\subseteq V(\frak{a})$ and $ Ext^{i}_{R}(\frac{R}{\frak{a}}, T)$ is finitely generated $R$-module for all $i\geq 0$.
\item[(ii)] $T$ is called a minimax $R$-module if there is a finitely generated submodule $L$ such that $T/L$ is Artinian $R$-module.
\item[(iii)] we say that $T$ is weakly laskerian $R$-module if the set of associated primes of any quotient module of $T$ is finite.
\item[(iv)] $T$ is $\frak{a}$-weakly cofinite $R$-module if $\hbox{Supp}M \subseteq V(\frak{a})$ and $\hbox{Ext}^{i}_{R}(\frac{R}{\frak{a}}, T)$ is weakly laskerian $R$-module for all $i\geq 0$.
 \end{itemize}
{\bf Remark 3.2}. We recall that the $N-\hbox{depth}$ of $L$, denoted by $\hbox{depth}(L, N)$, is defined as the length of any maximal $N$-sequence contained in $(0:_{R}L)$. Then $\hbox{depth}(L, N)$ is equal to the least integer $r$ such that $\hbox{Ext}^{r}_{R}(L, N)\neq 0$.\\
In addition, we denote by $\dim\hbox{Supp}H^{i}_{\frak{a}}(M, N)$ the maximum of numbers $\dim(R/\frak{p})$, where $\frak{p}$ runs over the Support of $H^{i}_{\frak{a}}(M, N)$.\\
 {\bf Theorem 3.3}. Let $t$ and $k$ be non-negative integers. If $\dim\hbox{Supp}H^{i}_{\frak{a}}(M, N)\leq k$ for all $i<t$, and put $\frak{a}+\hbox{Ann}(M)=I$, then $\dim\hbox{Supp}H^{i}_{\frak{a}}(M, N/\Gamma_{I}(N))\leq k.$\\
{\bf Proof}. From the short exact sequence $0\longrightarrow \Gamma_{I}(N)\longrightarrow N \longrightarrow N/\Gamma_{I}(N) \longrightarrow 0$, we get the long exact sequence $$H^{i}_{\frak{a}}(M, \Gamma_{I}(N))\longrightarrow H^{i}_{\frak{a}}(M, N) \longrightarrow H^{i}_{\frak{a}} (M,N/\Gamma_{I}(N)) \longrightarrow H^{i+1}_{\frak{a}}(M, \Gamma_{I}(N))$$ for all $i$. We assume that there exists an integer $i<t$ and $\frak{p}\in \hbox{Supp}H^{i}_{\frak{a}}(M, N/\Gamma_{I}(N))$ such that $\dim(R/\frak{p})>k$ and $\frak{p}\notin \hbox{Supp}H^{j}_{\frak{a}}(M, N/\Gamma_{I}(N))$ for all $j<i$. Thus by the long exact sequence as above, in conjunction with the fact that $H^{i}_{\frak{a}}(M, \Gamma_{I}(N))\cong \hbox{Ext}^{i}_{R}(M, \Gamma_{I}(N))$, we obtain the following exact sequence $$\hbox{Ext}^{i}_{R}(M, \Gamma_{I}(N))_\frak{p}\longrightarrow H^{i}_{\frak{a}}(M, N)_\frak{p} \longrightarrow H^{i}_{\frak{a}} (M,N/\Gamma_{I}(N)) _\frak{p}\longrightarrow \hbox{Ext}^{i+1}_{R}(M, \Gamma_{I}(N))_\frak{p}. \;\;\;\;(2)$$ Note that $H^{j}_{\frak{a}}(M, N)_\frak{p}=0$ for all $j\leq i$, while $ H^{j}_{\frak{a}} (M,N/\Gamma_{I}(N)) _\frak{p}=0$ for all $j<i$, and  $H^{i}_{\frak{a}} (M,N/\Gamma_{I}(N)) _\frak{p}\neq 0$. So, by the exact sequence (2), we have $\hbox{Ext}^{j}_{R}(M, \Gamma_{I}(N))_\frak{p}=0$ for all $j\leq i$, and $\hbox{Ext}^{i+1}_{R}(M, \Gamma_{I}(N))_\frak{p}\neq 0$. It implies that $\Gamma_{I}(N)_\frak{p}\neq 0$ and $\hbox{depth}(M_\frak{p}, \Gamma_{I}( N)_\frak{p})=i+1\geq 1$. Hence $\hbox{Ann}(M)_\frak{p}\nsubseteq \frak{q}R_\frak{p}$ for all $\frak{q}R_\frak{p}\in \hbox{Ass}_{R_{\frak{p}}}(\Gamma_{I}(N))_\frak{p}$. This contradicts with the fact that $\hbox{Ass}_{R_{\frak{p}}}(\Gamma_I(N))_\frak{p}=\hbox{Ass}_{R_{\frak{p}}}(N_\frak{p})\bigcap V(I)_{\frak{p}}$.\\
{\bf Lemma 3.4}. Let $T$ be a representable $R$-module and let $\hbox{Hom}(\frac{R}{\frak{a}}, T)$ is of finite length. Then $V(\frak{a})\cap \hbox{Att}(T)\subseteq \hbox{Max}(R)$.\\
{\bf Proof}. Let $\frak{p}\in V(\frak{a})\cap \hbox{Att}(T)$, and let $T=s_{1}+s_2+\cdots +s_{n}$, where $s_{r}$ is $\frak{p}$-secondary. It follows that there is a integer $k$ such that $\frak{p}^{k}s_{r}=0$, and so $\frak{a}^{k}s_{r}=0$. As $\hbox{Hom}(\frac{R}{\frak{a}}, T)$ is of finite length, $\hbox{Hom}(\frac{R}{\frak{a}^{k}}, T)$ is of finite length. Since $s_{1}\subseteq \hbox{Hom}(\frac{R}{\frak{a}^{k}}, T)$, it follows that $s_{1}$ has finite length. Therefore $\frak{p}\in\hbox{Max}(R)$, as required.\\
{\bf Lemma 3.5}. Let $T$ be a representable $R$-module, and $\frak{q}\in\hbox{Spec}(R)$. Let $x\notin {\bigcup}^{n}_{i=1}\frak{p}_{i}$, where $\frak{p}_{i}\in\hbox{Att}(T)-\{\frak{q}\}$. Then $\hbox{Att}(T/xT)\subseteq \{\frak{q}\}$.\\
{\bf Proof}. Let $s_{1}+\cdots +s_{n}+s$ be a minimal secondary representation of $T$, where $s$ is $\frak{q}$-secondary and $s_{i}$ is $\frak{p}_{i}$-secondary for all $i=1, \cdots, n$. Since $x\notin\bigcup^{n}_{i=1}\frak{p}_{i}$ for all $1 \leq i \leq n$, so $xs_{i}=s_{i}$. Note that $(T/xT)=s/(s\cap (s_{1}+\cdots+s_{n}+xs))$. Therefore $\hbox{Att}(T/xT) \subseteq \hbox{Att}s=\{\frak{q}\}$. On other hand, if $\frak{q}\notin\hbox{Att}(T)$, then $\hbox{Att}(T/xT)=\hbox{Att}(0)=\emptyset$, our claim is clear.\\
{\bf Theorem 3.6}. Let $M$, $N$ be two finitely generated $R$-module and let $t$ be non-negative integer such that $\hbox{Supp}H^{i}_{\frak{a}}(M, N)$  is finite set for all $i<t$. Then $H^{i}_{\frak{a}}(M, N)$ is minimax and $\frak{a}$-cofinite for all $i<t$ and $\hbox{Hom}(R/\frak{a}, H^{t}_{\frak{a}}(M, N))$ is finitely generated.\\
{\bf Proof}. Since the number of prime ideals between two given ones in Noetherian ring is zero or infinite, $\dim\hbox{Supp}H^{i}_{\frak{a}}(M, N)\leq 1$ for all $i<t$. Also, as $H^{i}_{\frak{a}+\hbox{Ann}(M)}(M, N)\cong H^{i}_{\frak{a}}(M, N)$ for all $i$, we can assume that $\hbox{Ann}(M)\subseteq V(\frak{a})$. Now, we prove the claim by induction on $t\geq 0$. The case of $t=0$ is clear. If $t=1$ then it is that $H^{0}_{\frak{a}}(M, N)=\hbox{Hom}(M, \Gamma_{\frak{a}}(N))$ is minimax and $\frak{a}$-cofinite and $\hbox{Hom}(R/\frak{a}, H^{1}_{\frak{a}}(M, N))$ is finitely generated by [6, Theorem 2.7]. Assume that $t>1$ and the result holds true for the case $t-1$. From the short exact sequence $0\longrightarrow \Gamma_{\frak{a}}(N)\longrightarrow N \longrightarrow N/\Gamma_{\frak{a}}(N)\longrightarrow 0$, we get the long exact sequence
$$H^{i}_{\frak{a}}(M, \Gamma_{\frak{a}} (N)) \overset{f_i} \longrightarrow H^{i}_{\frak{a}}(M, N)\overset{g_i} \longrightarrow H^{i}_{\frak{a}}(M, N/\Gamma_{\frak{a}}(N)) \overset{h_i} \longrightarrow  H^{i+1}_{\frak{a}}(M, \Gamma_{\frak{a}}(N)).$$ For each $i\geq 0$ we split the above exact sequence into two the following exact sequences
\begin{multline*}
0\longrightarrow \hbox{im}f_i \longrightarrow H^{i}_{\frak{a}}(M, N) \longrightarrow \hbox{im}g_i \longrightarrow 0 \;\;\;\;\;\hbox{and}\\
0\longrightarrow \hbox{im}g_i \longrightarrow H^{i}_{\frak{a}}(M, N/\Gamma_{\frak{a}}(N)) \longrightarrow \hbox{im}h_i \longrightarrow 0.
\end{multline*}
Note that $\hbox{im}f_i$ and $\hbox{im}h_i$ are finitely generated for all $i\geq 0$. Then, for each $i<t$, we obtain that $H^{i}_{\frak{a}}(M, N)$ is $\frak{a}$-cofinte and minimax module if and only if so is $H^{i}_{\frak{a}}(M, N/\Gamma_{\frak{a}}(N))$. On the other hand, we get by Theorem 3.3 that $\dim\hbox{Supp}H^{i}_{\frak{a}}(M, N/\Gamma_{\frak{a}}(N))\leq 1$ for all $i<t$. Therefore, in order to prove the theorem for the case of $t>1$, we may assume that $\Gamma_{\frak{a}}(N)=0$. In addition, put $X=\bigcup^{t-1}_{i=0}\hbox{Supp}(H^{i}_{\frak{a}}(M, N))$, $S=\{\frak{p}\in X\mid \dim(R/\frak{p}=1\}$. Thus $S\subseteq \bigcup^{t-1}_{i=0}\hbox{Ass}(H^{i}_{\frak{a}}(M, N))$.
By using the inductive hypothesis the $R$-module $\hbox{Hom}(R/{\frak{a}}, H^{t-1}_{\frak{a}}(M, N))$ is finitely generated and $R$-module $H^{i}_{\frak{a}}(M, N)$ is $\frak{a}$-cofinite and minimax for all $i<t-1$. It implies that $\bigcup^{t-1}_{i=0}\hbox{Ass}(H^{i}_{\frak{a}}(M, N))$ is a finite set, and so $S$ is a finite set. Assume that $S=\{\frak{p}_1, \frak{p}_2,\cdots , \frak{p}_{n}\}$. Now, it is straightforward to see that $\hbox{Supp}_{R_{\frak{p}_k}}H^{i}_{\frak{a}}(M, N)_{\frak{p}_k}\subseteq \hbox{Max}(R_{\frak{p}_k})$ and $\hbox{Hom}(R/\frak{a}, H^{i}_{\frak{a}}(M, N))_{\frak{p}_k}$ is finitely generated for all $i<t$ and $k=1,\cdots,n$. It follows that $\hbox{Hom}(R/\frak{a}, H^{i}_{\frak{a}}(M, N))_{\frak{p}_k}$ is Artinian for all $i<t$. As $H^{i}_{\frak{a}}(M, N)_{\frak{p}_k}$ is $\frak{a}R_{\frak{p}_{k}}$-torsion, it yields from Melkerssons$^{,}$ theorem [11, Theorem 1.3]
$H^{i}_{\frak{a}}(M, N)_{\frak{p}_k}$ is Artinian for all $i<t$ and $k=1,\cdots , n$. Therefore $\bigcup^{t-1}_{i=0}\hbox{Att}H^{i}_{\frak{a}}(M, N)_{\frak{p}_k}$ is finite set. Now, we choose an element
$x\in \frak{a}$ such that $x\notin\bigcup_{\frak{p}\in\hbox{Ass}(N)}\frak{p}\bigcup_{i=0}^{t-1}\bigcup_{\frak{p}\in Y_{i}} \frak{p}$, where $Y_i= \bigcup_{k=1}^{n}\{\frak{q}\mid \frak{q}R_{\frak{p}_k}\in \hbox{Att}H^{i}_{\frak{a}}(M, N)_{\frak{p}_{k}}\}-V(\frak{a})$. It follows that $x$ is an $N$-sequence. Therefore, we may consider the exact sequence $0\longrightarrow N \overset{x}\longrightarrow N \longrightarrow  N/xN \longrightarrow 0$ to obtain the exact sequence$$H^{i}_{\frak{a}}(M, N)\overset{x}\longrightarrow H^{i}_{\frak{a}}(M, N) \longrightarrow H^{i}_{\frak{a}}(M, N/xN) \longrightarrow H^{i+1}_{\frak{a}}(M, N)$$
for all $i\geq 0$. It implies the following exact sequence $$0\longrightarrow H^{i}_{\frak{a}}(M, N)/xH^{i}_{\frak{a}}(M, N)\overset{\alpha_{i}}\longrightarrow H^{i}_{\frak{a}}(M, N/xN) \overset{\beta_{i}}\longrightarrow (0:_{H^{i+1}_{\frak{a}}(M, N)}x)\longrightarrow 0\;\;\;\;(3)$$ for all $i\geq 0$. By using the exact sequence (3) in conjunction with the hypothesis, yields the $\dim\hbox{Supp}(H^{i}_{\frak{a}}(M, N/xN))\leq 1$ for all $i<t-1$. So that, we get by the induction assumption that $H^{i}_{\frak{a}}(M, N/xN)$ is $\frak{a}$-cofinite and minimax for all $i<t-1$ and $\hbox{Hom}(\frac{R}{\frak{a}}, H^{t-1}_{\frak{a}}(M, N/xN))$ is finitely generated.\\ Moreover, since $\hbox{Hom}(\frac{R}{\frak{a}}, H^{i}_{\frak{a}}(M, N))_{\frak{p}_{k}}$ is of finite length for all $i<t$ and $k=1,\cdots,n$. Therefore, by Lemma 3.4 $V(\frak{a})_{\frak{p}_{k}} \bigcap \hbox{Att}H^{i}_{\frak{a}}(M, N)_{\frak{p}_{k}}\subseteq \{\frak{p}_{k}R_{\frak{p}_{k}}\}$. By the choice of $x$ and Lemma 3.5, we obtain $(\frac{H^{i}_{\frak{a}}(M, N)}{xH^{i}_{\frak{a}}(M, N)})_{\frak{p}_{k}}$ has finite length for all $i<t$ and all $k=1, \cdots, n$. Then, there exists a finitely generated submodule $L_{i_{k}}$ of $L_{i}=(H^{i}_{\frak{a}}(M, N)/xH^{i}_{\frak{a}}(M, N))$ such that $(L_{i_{k}})_{\frak{p}_{k}} = (L_{i})_{\frak{p}_{k}}$. Put $V_{i}= L_{i_{1}}+L_{i_{2}}+\cdots +L_{i_{n}}$. Then $V_i$ is a finitely generated submodule of $L_{i}$ and $\hbox{Supp}(\frac{L_{i}}{V_{i}})\subseteq X-\{\frak{p}_{1}, \cdots , \frak{p}_{n}\} \subseteq \hbox{Max}(R)$ for all $i<t$. Since $\hbox{Hom}(\frac{R}{\frak{a}}, H^{i}_{\frak{a}}(M, N/xN))$ is finitely generated for all $i<t$, it follows from exact sequence $(3)$ that $\hbox{Hom}(\frac{R}{\frak{a}}, L_{i})$ is finitely generated for all $i< t$. Consideration of the exact sequence $0\longrightarrow V_i\longrightarrow L_{i}\longrightarrow L_{i}/V_{i} \longrightarrow 0$, shows that $\hbox{Hom}(\frac{R}{\frak{a}}, \frac{L_{i}}{V_{i}})$ is a finitely generated $R$-module for all $i<t$. It is also Artinian module, because it is supported only at maximal ideals. As $\frac{L_{i}}{V_{i}}$ is $\frak{a}$-torsion, it yields from Melkersson$^{,}$ theorem [11, Theorem 1.3], that $R$-module $\frac{L_{i}}{V_{i}}$ is Artinian, and so $L_{i} $ is minimax for all $i<t$. Now, by [13, Proposition 4.3] in conjunction with the fact that $\hbox{Hom}(\frac{R}{\frak{a}}, L_{i})$ is finitely generated $R$-module for all $i<t$, to see that $L_{i}$ is $\frak{a}$-cofinite and minimax $R$-module for all $i<t$. Again, consideration of the exact sequence $(3)$ shows that $(0:_{H^{t-1}_{\frak{a}}(M, N)}x)$ is $\frak{a}$-cofinite and minimax and $\hbox{Hom}(\frac{R}{\frak{a}}, (0 :_{H^{t}_{\frak{a}}(M, N)}x)$ is finitely generated $R$-module. It implies that $H^{t-1}_{\frak{a}}(M, N)$ is $\frak{a}$-cofinite and minimax by [13, Proposition 4.3]. The following completes the proof:$$\hbox{Hom}(\frac{R}{\frak{a}}, H^{t}_{\frak{a}}(M, N)= \hbox{Hom}(\frac{R}{\frak{a}}\otimes \frac{R}{(x)}, H^{t}_{\frak{a}}(M, N))=(0:_{H^{t-1}_{\frak{a}}(M, N)}x).$$
{\bf Theorem 3.7}. Let $(R, \frak{m})$ be local ring and let $t$ be a non-negative integers such that $\dim\hbox{Supp}H^{i}_{\frak{a}}(M, N)\leq 2$ for all $i<t$. Then $H^{i}_{\frak{a}}(M, N)$ is $\frak{a}$-weakly cofinite for $j\geq 0$ and $i<t$.\\
{\bf Proof}. Since $\hbox{Supp}H^{i}_{\frak{a}}(M, N)\subseteq V(\frak{a})$, it is enough to show that $\hbox{Ext}^{j}_{R}(\frac{R}{\frak{a}}, H^{i}_{\frak{a}}(M, N))$ is weakly laskerian for all $j \geq 0$ and $i< t$. Using [3, Theorem 4.3.2] and [9, Ext 7.7], without losing generality we may assume that $R$ is complete. Now, suppose, contrary to our claim, that there is fixed integers $i<t$ and $j\geq 0$ such that $\hbox{Ext}^{j}_{R}(\frac{R}{\frak{a}}, H^{i}_{\frak{a}}(M, N))$ is not weakly laskerian. In view of definition, let $T^{\prime}$ be a submodule of $T=\hbox{Ext}^{j}_{R}(\frac{R}{\frak{a}}, H^{i}_{\frak{a}}(M, N))$ and $\hbox{Ass}(T/T^{\prime})$ is infinite set. Then, the set $ \hbox{Ass}(\frac{T}{T^{\prime}})-\{\frak{m}\}$ is not empty. By [10, Lemma 3.2], there exists $x\in\frak{m}$ such that $x\notin\bigcup_{\frak{p}\in A}\frak{p}$. Let $S=\{x^{k}\mid 0\leq k \in \BZ\}$. Then $S\cap\frak{m}\neq \emptyset$. So, $\dim\hbox{Supp}S^{-1}(H^{i}_{\frak{a}}(M, N))\leq 1$. It follows from Theorem 3.5, $\hbox{Ext}^{j}_{S^{-1}R}(\frac{S^{-1}R}{S^{-1}\frak{a}}, S^{-1}(H^{i}_{\frak{a}}(M, N)))$ is  finitely generated. Therefore,  $\hbox{Ass}\frac{S^{-1}T}{S^{-1}T^{\prime}}$ is finite set, and so by [9, Theorem 6.2] $\hbox{Ass}S^{-1}(T/T^{\prime})=\{S^{-1}\frak{p}\mid \frak{p}\in\hbox{Ass}T/T^{\prime}, \frak{p}\cap S=\varnothing \}$, which is a contradiction.
\section*{4. Artinianness of composed graded local cohomology }
The concept the tameness is the most fundamental concept related to the asymptotic behaviour of cohomology modules. A graded $R$-module $T=\bigoplus T_{n}$ is said to be tame or asymptotic gap-free if either $T_{n}=0$ for all $n\ll 0$ else $T_{n}\neq 0$ for all $n\ll 0$. It is known that any graded Artinian $R$-module is tame. In this section, we keep the notation and hypotheses introduced in the introduction and we study the Artinianness of graded module $H^{j}_{\frak{b}_{0}}(H^{i}_{\frak{a}}(M, N))$ for some $i\geq 0$ and $j\geq 0$.\\
{\bf Theorem 4.1}. Let $T$ be an $\frak{a}$-torsion minimax graded $R$-module and $\sqrt{\frak{a}_{0}+\frak{b}_{0}}=\frak{m}_{0}$. Then $H^{i}_{\frak{b}_{0}}(T)$ and $Tor^{R_{0}}_{i}(\frac{R_{0}}{\frak{b}_{0}}, T)$ are Artinian for all $i\geq 0$.\\
{\bf Proof}. As $T$ is graded minimax $R$-module, there is an exact sequence of graded $R$-modules $$0\longrightarrow T^{\prime} \longrightarrow T \longrightarrow T/T^{\prime} \longrightarrow 0\;\;\;\;\;\;(4)$$ such that $T^{\prime}$ is graded finitely generated and $T/T^{\prime}$ is graded Artinian. The application of local cohomology with respect to $\frak{b}_{0}$ to it leads to an exact sequence of graded $R$-modules $0 \longrightarrow \Gamma_{\frak{b}_{0}}(T^{\prime})\longrightarrow \Gamma_{\frak{b}_{0}}(T) \longrightarrow \Gamma_{\frak{b}_{0}}(T/T^{\prime}) \longrightarrow H^{1}_{\frak{b}_{0}}(T^{\prime}) \longrightarrow H^{1}_{\frak{b}_{0}}(T)\longrightarrow 0$ and the isomorphisms $H^{i}_{\frak{b}_{0}}(T^{\prime})\cong H^{i}_{\frak{b}_{0}}(T)$ for all $i\geq 2$. We note that $\Gamma_{\frak{b}_{0}}(T/T^{\prime})=T/T^{\prime}$ is Artinian. Since $\Gamma_{\frak{a}}(T)=T$ , the submodule $T^{\prime}$ is $\frak{a}$-torsion. Then $H^{i}_{\frak{b}_{0}}(T^{\prime})\cong H^{i}_{\frak{m}}(T^{\prime})$ for all $i\geq 0$. This proves $H^{i}_{\frak{b}_{0}}(T)$ is Artinian for all $i\geq 0$. For the second claim, if we apply the functor $Tor^{R_{0}}_{i}(\frac{R_{0}}{\frak{b}_{0}}, \;\;\;\;)$ to the short exact sequence (4), we have the following exact sequence of graded $R$-modules
$$Tor^{R_{0}}_{i+1}(\frac{R_{0}}{\frak{b}_{0}}, T/T^{\prime})\longrightarrow Tor^{R_{0}}_{i}(\frac{R_{0}}{\frak{b}_{0}}, T^{\prime})\longrightarrow Tor^{R_{0}}_{i}(\frac{R_{0}}{\frak{b}_{0}}, T) \longrightarrow Tor^{R_{0}}_{i}(\frac{R_{0}}{\frak{b}_{0}}, T/T^{\prime}).\;\; (5)$$According to [6, Lemma 2.1], the graded module $Tor^{R_{0}}_{i}(\frac{R_{0}}{\frak{b}_{0}}, T/T^{\prime})$ is Artinian for each $i$. Since $Tor^{R_{0}}_{i}(\frac{R_{0}}{\frak{b}_{0}}, T^{\prime})$ is an $\frak{a}$-torsion finitely generated graded $R$-module,  $\frak{m}^{k_{i}}Tor^{R_{0}}_{i}(\frac{R_{0}}{\frak{b}_{0}}, T^{\prime})=0$ for some $k_{i}\in \BN$. It follows that $Tor^{R_{0}}_{i}(\frac{R_{0}}{\frak{b}_{0}}, T^{\prime})$ is Artinian for all $i\geq 0$. In view of exact sequence (4), $Tor^{R_{0}}_{i}(\frac{R_{0}}{\frak{b}_{0}}, T)$ is Artinian.\\
{\bf Theorem 4.2}. Let $T$ be an $\frak{a}$-torsion and $\frak{a}$-cofinite graded $R$-module and $\sqrt{\frak{a}_{0}+\frak{b}_{0}}=\frak{m}_{0}$. Then $H^{i}_{\frak{b}_{0}}(T)$ is Artinian and $\frak{a}$-cofinite $R$-module for all $i\geq 0$.\\
{\bf Proof}. We proceed by induction on $i$. If $i=0$, then $\Gamma_{\frak{b}_{0}} (T)=\Gamma_{\frak{m}}(T)$
is Artinian and $\frak{a}$-cofinite by [12, Corollary 1.8]. Let $i=1$. As $T$ is $\frak{a}$-torsion, there is a short exact sequence of $\frak{a}$-torsion graded modules $0\longrightarrow T\longrightarrow E\overset{d}\longrightarrow C \longrightarrow 0$ such that $E$ is injective and $C$ is $\frak{a}$-cofinite as well. Application of the functor $\Gamma_{\frak{b}_{0}}(\;\;\;)$ induces the following exact sequence
$$0\longrightarrow \Gamma_{\frak{b}_{0}}(T)\longrightarrow \Gamma_{\frak{b}_{0}}(E)\overset{\Gamma_{\frak{b}_{0}}(d)}\longrightarrow \Gamma_{\frak{b}_{0}}(C) \longrightarrow H^{1}_{\frak{b}_{0}}(T)\longrightarrow o\;\;.$$ Consider $X=\hbox{im}(\Gamma_{\frak{b}_{0}}(d))$. In view of the case $i=0$, we deduce that $\Gamma_{\frak{b}_{0}}(T)$ and $\Gamma_{\frak{b}_{0}}(C)$ are Artinian and $\frak{a}$-cofinite. So $\hbox{Hom}(\frac{R}{\frak{a}}, X)$ is finitely generated and hence the isomorphism $\hbox{Ext}^{i}_{R}(\frac{R}{\frak{a}}, X)\cong \hbox{Ext}^{i+1}_{R}(\frac{R}{\frak{a}}, \Gamma_{\frak{b}_{0}}(T))$, for each $i\geq 1$, implies that $X$ is Artinian and $\frak{a}$-cofinite. Now, the cofiniteness and Artinianness of $X$ and $\Gamma_{\frak{b}_{0}}(C)$ imply that $H^{1}_{\frak{b}_{0}}(T)$ is Artinian and $\frak{a}$-cofinite. On the other hand, for each $i>1$, there is an isomorphism $H^{i-1}_{\frak{b}_{0}}(C)\cong H^{i}_{\frak{b}_{0}}(T)$. Therefore the induction completes the proof.\\
{\bf Theorem 4.3}. Let $t$ be a non-negative integer such that $\hbox{Supp} H^{i}_{\frak{a}}(N)$ is finite set for all $i<t$ and $\sqrt{\frak{a}_{0}+\frak{b}_{0}}=\frak{m}_{0}$. Then $H^{j}_{\frak{b}_{0}}(H^{i}_{\frak{a}}(M, N))$ is Artinian and $\frak{a}$-cofinite for all $i<t$ and $j\geq 0$. In addition,  $H^{j}_{\frak{b}_{0}}(H^{t}_{\frak{a}}(M, N))$ is Artinian  for all $j=0, 1$.\\
{\bf Proof}. In view of Theorem 2.1, there is a Grothendieck$^{,}$s spectral sequence $$E_{2}^{p,q}=\hbox{Ext}^{p}_{R}(M, H^{q}_{\frak{a}}(N))\underset{p}\Longrightarrow H^{p+q}_{\frak{a}}(M, N).$$
 Using Theorem 3.4 $ H^{i}_{\frak{a}}(N)$ is minimax and $\frak{a}$-cofinite for all $i<t$. It follows from Theorem 2.3, $H^{i}_{\frak{a}}(M, N)$ is minimax and $\frak{a}$-cofinite for all $i<t$. Then, from Theorem 4.2, $H^{j}_{\frak{b}_{0}}(H^{i}_{\frak{a}}(M, N))$ is Artinian and $\frak{a}$-cofinite for all $i<t$ and all $j\geq 0$. On the other, using the spectral sequence $$E_{2}^{p,q}=H_{\frak{b}_{0}}^{p}(H^{q}_{\frak{a}}(M, N)\underset{p} \Longrightarrow H^{p+q}_{\frak{m}}(M, N)$$ in conjunction Theore 2.4 and the fact that $H^{i}_{\frak{m}}(M, N)$ is Artinian for all $i\geq 0$, the result follows.\\
{\bf Corollary 4.4}. Let $t$ be a non-negative integer such that $\dim(\hbox{Supp} H^{i}_{\frak{a}}(M, N))$ is finite set for all $i<t$ and $\sqrt{\frak{a}_{0}+\frak{b}_{0}}=\frak{m}_{0}$. Then $H^{j}_{\frak{b}_{0}}(H^{i}_{\frak{a}}(M, N))$ is Artinian and $\frak{a}$-cofinite for all $i<t$ and $j\geq 0$. In addition,  $H^{j}_{\frak{b}_{0}}(H^{t}_{\frak{a}}(M, N))$ is Artinian  for all $j=0, 1$.\\
{\bf Proof}. It follows immediately by using Theorem 4.3.\\
{\bf Definition 4.5}. A sequence $x_{1},\cdots\;,x_{n}$ of elements of $\frak{a}$ is said to be a generalized regular sequence of $T$ if $x_{i}\notin \frak{p}$ for all $\frak{p}\in \hbox{Ass}_{R}(T/(x_{1},\cdots\;,x_{i-1})T)$ satisfying $\dim(R/\frak(p)> 1$ for all $i=1,2,\cdots\;,n$. It is clear that if $\dim(T/\frak{a}T >1$, each generalized regular sequence of $T$ in $\frak{a}$ has finite length. The length of a maximal regular sequence of $T$ in $\frak{a}$ is denoted by $\hbox{gdepth}(\frak{a}, T)$. [14]\\
Also, $\hbox{gdepth}(M/\frak{a}M, N)$ is a non-negative integer and is equal to the length of any maximal generalized regular $N$-sequence in $\frak{a}+(0: M)$.\\
{\bf Theorem 4.6}. Let $\hbox{gdepth}(M/\frak{a}M, N)=cd_{\frak{a}}(M, N)$, where $cd_{\frak{a}}(M, N)=\hbox{sup}\{i\in\BN_{0}\mid H^{i}_{\frak{a}}(M, N)\neq 0\}$. Then $H^{j}_{\frak{b}_{0}}(H^{i}_{\frak{a}}(M, N))$ is Artinian and cofinite for all $i\geq 0$ and $j\geq 0$. In particular $H^{j}_{\frak{b}_{0}}(H^{i}_{\frak{a}}(M, N))$ is Artinian and $\frak{a}$-cofinite for all $i\geq 0$ and $j\geq 0$, if $N$ is $\frak{a}$-cofinite.\\
{\bf Proof}. Using [16, Theorem 3.2], $\hbox{Supp} H^{i}_{\frak{a}}(M, N)$ is finite set for all $i<\hbox{gdepth}(M/\frak{a}M, N)$. It follows from Corollary 4.4 that $R$-module $H^{j}_{\frak{b}_{0}}(H^{i}_{\frak{a}}(M, N))$ is Artinian and $\frak{a}$-cofinite for all $i< \hbox{gdepth}(M/\frak{a}M, N)$ and all $j\geq 0$. If $i> cd_{\frak{a}}(M, N)$, then, in view of the definition of $cd_{\frak{a}}(M, N)$, $H^{j}_{\frak{b}_{0}}(H^{i}_{\frak{a}}(M, N))=0$. Thus we consider the case where $i=t=\hbox{gdepth}(M/\frak{a}M, N)$. To this end, consider the Grothendieck$^{,}$s spectral sequence $$E_{2}^{p,q}=H^{p}_{\frak{b}_{0}}(H^{q}_{\frak{a}}(M, N))\underset{p}\Longrightarrow \;\;H^{p+q}_{\frak{m}}(M, N).$$
 According to Theorem 2.1(ii), we have $H^{j}_{\frak{b}_{0}}(H^{t}_{\frak{a}}(M, N))\cong (H^{t+j}_{\frak{m}}(M, N))$. In view of [3, Theorem 7.1.3], $H^{i}_{\frak{m}}(N)$ is Artinian $R$-module and by [8, Lemma 2.7], $H^{i}_{\frak{m}}(N)$ is Artinian and $\frak{a}$-cofinite $R$-module if $N$ is $\frak{a}$-cofinite $R$-module. Therefore, by using Grothendieck$^{,}$s spectral sequence $$ E_{2}^{p,q}=\hbox{Ext}_{R}^{p}(N, H^{q}_{\frak{m}}(N))\underset{p}\Longrightarrow H_{\frak{m}}(M, N), $$the result follows from Theorem (2.3).

{\bf Acknowledgment.}
 I would like to thank professors Amir Mafi for their careful reading of the first draft and many helpful suggestions.


{\bf Fatemeh Dehghani-Zadeh }\\
Department of Mathematics ,
Islamic Azad University \\
Yazd Branch ,
Yazd, Iran \\
e-mail: f.dehghanizadeh@yahoo.com\\
\mbox{} e-mail: fdzadeh@gmail.com
\end{document}